\documentstyle[amscd,amssymb,verbatim,12pt]{amsart}
\input xypic
\xyoption{all}
\input epsf

\def\hcorrection#1{\advance\hoffset by #1 }
\def\vcorrection#1{\advance\voffset by #1 }

\vcorrection{-.5in}
\hcorrection{-.5in}

\newcommand{\B}[1]{{\bold#1}} 
\newcommand{\C}[1]{{\cal#1}} 

\theoremstyle{plain}
\newtheorem{theorem}{Theorem}[section]
\newtheorem{cor}{Corollary}[section]

\newtheorem{lem}{Lemma}[section]
\newtheorem{prop}{Proposition}[section]

\theoremstyle{definition}
\newtheorem{defin}{Definition}[section]

\theoremstyle{definition}
\newtheorem{example}{Example}[section]

\theoremstyle{remark}

\newtheorem{notation}{Notation}[section]

\numberwithin{equation}{section}

\begin{document}
\pagestyle{plain}
\addtolength{\footskip}{.3in}

\title{Hochschild DGLAs and Torsion Algebras}
\author{Lucian Ionescu}
\address{Mathematics Department\\Kansas State University\\
             Manhattan, Kansas 66502}
\email{luciani@@math.ksu.edu}
\thanks{I thank Jim Stasheff for helpful remarks,
and for pointing to the related work of \cite{DMM,DM}.}
\date{Octomber 1999}

\begin{abstract}
The associator of a non-associative algebra is the curvature of the
Hochschild quasi-complex. 
The relationship ``curvature-associator'' is investigated.
Based on this generic example,
we extend the geometric language of vector fields
to a purely algebraic setting,
similar to the context of Gerstenhaber algebras.
We interprete the elements of a non-associative algebra with
a Lie bracket as ``vector fields'' and the multiplication as a connection.
Conditions for the existance of an ``algebra of functions''
having as algebra of derivations the original non-associative algebra,
are investigated.
\end{abstract}

\maketitle

\tableofcontents

\newcommand{\one}{{\bf 1}}
\newcommand{\tensor}[1]{\underset{#1}{\otimes}}
\newcommand{\To}{\longrightarrow}


\section{Introduction}\label{S:I}
The associator $\alpha$ of a non-associative algebra $(A,\mu)$
is the curvature of the Hochschild quasi-complex.
It has properties analog to the properties of the curvature of
a linear connection, suggesting the interpretation of its elements as
formal vector fields.
Gerstenhaber algebras motivate this interpretation, as a
non-commutative analog of Poisson algebras.

In this article we investigate the possibility of ''reconstructing''
an algebra of ``functions'' representing a ''non-commutative space'',
such that the original algebra is the algebra of
its derivations (''vector fields'').
A differential calculus based on derivations is only sketched.
Such a calculus was developed in detail in \cite{DMM,DM}.

\vspace{.1in}
Classical differential geometry is built on the notion of {\em space}:
differential manifolds. A rough ``hierarchy'' is: space, functions,
vector fields and differential forms, connections etc.

Algebraic geometry starts at the ``second level'' (functions) by
considering an arbitrary commutative algebra and then constructing
the ``first  level'', the substitute for a space: its spectrum.
A {\em space} (affine variety) is roughly a pair consisting of a
topological space and its algebra of functions.

A natural question arises: ``What can be derived starting from the
``third level'' - vector fields - and to what extent is it profitable?''

We will be interested in geometry with connections, so we will investigate
the interpretation of the associator of a non-associative algebra
as a curvature (section \ref{SS:motiv}).
We recall in section \ref{S:hc}
some well know facts about the Lie algebra structure on the
Hochschld quasi-complex ($d^2\ne 0$) of a non-associative algebra.
Motivated by this generic example, we consider a Lie algebra,
with elements thought of as vector fields,
with a multiplication (not necessarily associative)
thought of as a linear connection,
and we define in section \ref{S:TA} the algebra of functions
($2^{\text{nd}}$ level).
Conditions are derived for the existence of a natural
``algebra of functions'', having as derivations
the original non-associative algebra (theorem \ref{th1}).
Section \ref{S:Ex} considers examples of torsion algebras.
The main example is the pre-Lie algebra of Hochschild cochains with
Gerstenhaber composition (theorem \ref{th2}).


As motivation for the emphasis on vector fields we mention two
sources.
Physical understanding evolved from considering
phase spaces (Poisson manifolds) rather than configuration spaces.
Moreover, the actual goal is to model the {\em space of evolutions}
of a system.

The second motivation is the correspondence between Poisson-Lie group
structures and Lie bialgebras $g$.
After quantization of $Ug$
one has deformed differential operators and,
in some sense deformed vector fields.
It would be convenient to have a procedure to recover an
algebra of ``functions'' (possibly non-commutative).

In deformation quantization of Poisson manifolds, one keeps the classical
observables and deforms the ``laws of mechanics'' to account for the
Heisenberg bracket.
We consider this approach as slightly conservative
and consider that the basic level for quantum physics: states - vectors,
evolution - operators, does not need an actual (configuration) space.
Our ``functions'' are naturally operators
on the given algebra of ``vector fields''.

%
%
\subsection{Associator - curvature / monoidal structure}\label{SS:motiv}
Another argument for a geometric interpretation of an algebra $(A,\mu)$
is as follows.
The associator $\alpha(x,y,z)=(xy)z-x(yz)$ of the
nonassociative algebra is formally a ``curvature''
of the left regular quasi-representation
$L:(A,\mu)\to (End_k(A),\circ)$:
\begin{equation}
L(xy) \stackrel{\sigma_{x,y}}{\To} L(x)L(y)
\end{equation}
$$\sigma(x,y)(z)=(L(x)L(y)-L(xy))(z)=-\alpha(x,y,z)$$
A {\em quasi-representation} is a morphism in the underlying category,
not necessarely preserving the additional structure
(e.g. commuting with the group operation).
In the situation at hand, the above map $L$ is assumed to be
only $k$-linear.
At the infinitesimal (Lie algebra) level, the same map $L$,
interpreted as a ``quasi-Lie'' representation
$L:(A,[,]_A)\to (End_k(A),[,])$,
defines a formal curvature $K$:
\begin{equation}
L([x,y]) \stackrel{K_{x,y}}{\To} [L(x),L(y)],
\qquad\qquad K(x,y)=[L(x),L(y)]-L([x,y])
\end{equation}
Moreover, in the Hochschild quasi-complex $(C^\bullet(A), d_\mu)$,
the differential $d_\mu$ has properties analog to a covariant derivative,
e.g. $d_\mu^2f=[\alpha,f]$.
Also $\alpha=\frac12[\mu,\mu]$ (the ``curvature'') is closed, 
i.e. a ``Bianchi identity'' holds. (see \ref{P:assoc}).

This approach is only a tentative to a {\em partial model} of
a ``local'' differential geometry.
The ``global'' point of view should consider the non-linear (multiplicative)
interpretation of the failure of a quasi-representation (or underlying morphism)
to preserve the structure as a non-strict monoidal structure.
This ``monoidal'' interpretation of the associator comes from looking at
quasi-representations, and more general at functions defined on groups
$s:G\to N$ as defining a 2-cocycle
$f:G\times G\to N, f(a\otimes b)=s(a)s(b) s(ab)^{-1}$
(``factor set'' corresponding to a given action of $G$ on $N$).
If interpreting $G$ and $N$ as monoidal groupoids \cite{I}, then $f$ is a
monoidal structure $f_{a,b}:s(a\otimes b)\to s(a)\otimes s(b)$
of the monoidal functor $s$.

\section{Hochschild Quasi-complex}\label{S:hc}
We begin by recalling the Hochschild cohomology \cite{Ge,GS} in the more
general case of a possibly non-associative algebra.

The obstruction to the cohomological study of a non-associative
algebras, $d^2\ne 0$, is interpreted as a curvature.
The use of the language of differential geometry is considered.

Throughout this section, $A$ will denote a module over a
commutative ring $R$.
We will assume that 2 and 3 do not annihilate non-zero elements in A.

\subsection{Pre-Lie Algebra of a Module}\label{SS:23}
Consider
$$C^{p,q}(A,A)
 =\{f:A^{\otimes^p}\to A^{\otimes^q} \vert f \ \text{R-linear}\}
 \qquad p,q\ge 0$$
and
$C^{\bullet,\bullet}(A,A)=\underset{p,q\ge 0}{\oplus} C^{p,q}$
with total degree $\partial eg(f^{p,q})=p-q$.
$C^{0,q}$ is identified with $A^{\otimes^q}$.
We will be interested in the first column
$C^\bullet (A)=\underset{p\in\B{N}}{\oplus} C^{p,1}$.
The grading induced by the total degree is
$C^{p-1}(A)=C^{p,1}(A,A)$ with $p\ge 0$.

We recall briefly the comp operation and the Lie algebra structure
it defines on the graded module of Hochschild cochains \cite{Ge,KPS}.

If $f^p\in C^p(A)$ and $g^q\in C^q(A)$, define the composition into the
$i^{th}$ place, where $i=1,\dots,p+1$:
$$f^p\circ_i g^q (a_1,\dots,a_{p+q-1})=f^p(a_1,\dots,a_{i-1},
g^q(a_i,\dots,a_{i+q-1}),a_{i+q},\dots,a_{p+q-1})$$
and the {\em comp} operation:
$$ f^p\circ g^q= \sum_{i=1}^{p+1} (-1)^{(i-1)q} f^p\circ_i g^q \in C^{p+q}$$
It is assumed that the composition is zero whenever $p=-1$.
Note that the (non-associative)  composition respects the grading.
Denote by $\alpha(f,g,h)=(f\circ g)\circ h - f\circ(g\circ h)$ the
{\em associator} of $\circ$, as a ``measure'' of non-associativity
of the comp oparation. 

The graded commutator is defined by:
$$[f^p,g^q]=f^p\circ g^q -(-1)^{pq}g^q\circ f^p$$
It is graded commutative:
$$[f^p,g^q]=-(-1)^{pq}[g^q,f^p]$$
and the graded Jacobi identity holds:
\begin{align}
(-1)^{FH}[f,[g,h]]+(-1)^{GF}[g,[h,f]]+(-1)^{HG}[h,[f,g]]=0
\end{align}
where $F,G,H$ denote the degrees of $f,g,h$ respectively. It is
equivalent to $ad$ being a representation of graded Lie algebras:
\begin{align} \label{E:derivation}
[f,[g,h]]=[[f,g],h]+(-1)^{pq}[g,[f,h]]\notag \\
ad_f([g,h])=[ad_f(g),h]+(-1)^{pq}[g,ad_f(h)]
\end{align}
%
\begin{notation}
We denote by $f\underset{(i,j)}{\circ} (g,h)$ the
simultaneous insertion of two functions $g$ and $h$ in the $i^{th}$ and
$j^{th}$ arguments of $f$ respectively.
\end{notation}
\begin{lem}\label{L:comp1}
If $f,g,h$ have degrees $p,q,r$, then:\\
(i) $f\circ_i(g\circ_j h)=(f\circ_i g)\circ_{i-1+j} h,
 \qquad 1\le i\le p+1,\quad 1\le j\le q+1$\\
(ii) $(f\circ g)\circ h - f\circ(g\circ h)=
  \sum_{i\ne j} \epsilon(i,j)
  (-1)^{(i-1)q+(j-1)r} f\underset{(i,j)}{\circ} (g,h)$, where
  $\epsilon(i,j)=1$ if $1\le j<i\le p+1$ and equals $(-1)^{qr}$ if
$1\le i<j\le p+1$.\\
(iii) $\alpha(f,g,h)=(-1)^{qr} \alpha(f,h,g)$
\end{lem}
\begin{pf}
(i) and (ii) follow from a straightforward inspection of trees
and signs. The key is that the only trees which survive, build out of $f,g,h$ in the
associator $\alpha$, are of the type $f\underset{(i,j)}{\circ} (g,h)$. The
``supercommutativity'' sign $(-1)^{qr}$ appears when $i$ passes over $j$ and
the order of insertion ($g$ before $h$) changes.
\end{pf}
\begin{notation}
Let $\mu\in C^1(A)$ and
$\mu=\mu_-+\mu_+$ the natural decomposition, with
$\mu_-(a,b)=\mu(a,b)-(-1)^{pq}\mu(b,a)$ and
$\mu_+(a,b)=\mu(a,b)+(-1)^{pq}\mu(b,a)$ the graded {\em skew}
and {\em symmetric} part of $\mu$. Alternatively $\mu_-$ will be
denoted as $[,]_\mu$ or just $[,]$ if no confusion is expected.
The corresponding associator will be denoted as $\alpha_\mu$.
\end{notation}
\begin{defin}
A (possibly non-associative) algebra $(A,\mu)$ is called a
{\em pre-Lie algebra} if $\mu_-$ is a Lie bracket.
\end{defin}
\begin{lem}
Let $(A,\mu)$ be an algebra and $\alpha$ its associator.\\
(i) $Alt(\alpha_{\mu_+})=0$\\
(ii) $Alt(\alpha_{\mu_-})=4 Alt(\alpha)$\\
(iii) $(A,\mu)$ is a pre-Lie algebra iff $Alt(\alpha)=0$.

If $A$ is graded, then a graded alternation $Alt$ is assumed.
\end{lem}
\begin{pf}
A direct computation.
\end{pf}
>From lemma \ref{L:comp1} and the above lemma it imediatly follows
the well known fact that the comp operation on Hochschild cochains
defines a Lie bracket.
\begin{cor}
$(C^\bullet,\circ)$ is a (graded) pre-Lie algebra.
\end{cor}
\begin{pf}
Since $\alpha(f,g,h)=(-1)^{qr}\alpha(f,h,g)$, we have:
$$Alt(\alpha)(f,g,h)=\sum_{cycl}\epsilon(f,g,h) (\alpha(f,g,h)-
(-1)^{qr}\alpha(f,h,g))=0$$
where $\epsilon(f,g,h)$, in the graded case, is not necessarily 1.
For example $\epsilon(g,h,f)=(-1)^{(q+r)p}$.
\end{pf}
%
\subsection{Quasi-complex of a Non-Associative Algebra}\label{S:hqc}
\begin{defin}
A {\em quasi-complex} is a sequence of objects and
morphisms in a category $\C{A}$
$$C^\bullet=\{\dots\to C^{-1}\to C^0\to C^1\to\dots\}.$$
The family of morphisms $d^\bullet$ is called a {\em quasi-differential},
for which $d^2$ may be non-zero. 
\end{defin}
Now assume that an element $\mu:A\otimes A\to A$ of degree one is fixed.
Thus $(A,\mu)$ is a (non-associative) R-algebra.
Define the quasi-differential as the adjoint $\mu$ action:
\begin{equation}
d_\mu(f^p)=[\mu,f]
\end{equation}
Then $(C^\bullet(A),d_\mu)$ is a {\em quasi-differential graded Lie algebra},
called the {\em Hochschild quasi-complex} corresponding to the algebra
$(A,\mu)$.

Note the difference of sign when compared with \cite{Ge,GS}:
$$d_{Ge}=-[f,\mu]=(-1)^p[\mu,f]=(-1)^p d_\mu f$$
As an example, for $p=1$ (with $d=d_\mu$ and $\mu(x,y)=xy$
if no confusion is possible):
\begin{align}\label{E:df}
df(x,y,z) &=\mu\circ f(x,y,z)+f\circ \mu (x,y,z)\notag\\
        &=\mu (f(x,y),z)-\mu (x,f(y,z))+f(\mu(x,y),z)-f(x,\mu (y,z))\notag\\
        &=-\{xf(y,z)-f(xy,z)+f(x,yz)-f(x,y)z\}
\end{align}
which is the usual Hochschild differential modulo the sign:
$$d_{Hoch}=(-1)^{\partial eg}\cdot d_\mu=d_{Ge}$$
Note that $d_\mu:C^p\to C^{p+1}$ has degree one,
$[ , ]:C^p\otimes C^q\to C^{p+q}$ is of degree zero and
$d_{Hoch}$ is not a graded derivation, since it does not satisfy
the ``Leibniz identity'' \ref{E:derivation}.
We state the following fact about graded Lie algebras, which is an imediate
consequence of the graded Jacobi identity.
\begin{lem}\label{L:1} 
Let $(g,[,])$ a graded Lie algebra over a commutative ring R.
Then, if $x$ is an even degree element $[x,x]=0$.
If $x$ is odd, $[x,[x,x]]=0$ and $ad_{[x,x]}=2 (ad_x)^2$.
\end{lem}
If the multiplication $\mu$ is associative, then (\ref{E:df}):
\begin{equation}\label{E:associator}
d\mu (x,y,z)=2\{(xy)z)-x(yz)\}=0
\end{equation}
and $[\mu,\mu]=0$. By the previous lemma 
$2[\mu,[\mu,f]]=d_{[\mu,\mu]}(f)=0$, and
thus $(C^\bullet(A),d_\mu)$ is a complex of R-modules.

We introduce the folowing definition:
\begin{defin}\label{D:assoc}
A quasi-differential graded Lie algebra $(C^\bullet, d)$ is
called {\em coboundary} if there is a degree zero element $I$ such that
the {\em multiplication} (also called {\em torsion})
$\mu=d I$ verifies the relation $d=ad_{\mu}$. 

Then $2\alpha=d\mu$ defines the {\em associator} $\alpha$
(or {\em curvature}).
\end{defin}
Since $d\mu=[\mu,\mu]=2\mu\circ\mu$
and 2 was assumed not to have right divisors, the associator is
well-defined by the above equation.

The motivation for the geometric terminology comes from a formal analogy
in the context of a derivation law in an $A$-module,
where $A$ is an $R$-algebra (see \cite{Ko})
or a linear connection $D$ on a vector bundle, where the torsion $T$ and
curvature $F$ of the total covariant derivative $d$ are defined as:
$$T=dI, \qquad dI(X,Y)=D_X I(Y) -D_Y I(X) - I([X,Y])$$
$$F(X,Y)=[D_X,D_Y]-D_{[X,Y]},\qquad d^2s=[F,s]$$
Here $I$ is the identity tensor.

The geometric interpretation will be considered in section \ref{S:TA}.
\begin{prop}\label{P:assoc}
Let $(A,\mu)$ be an R-algebra, possibly non-associative.
Then $(C^\bullet(A),d_\mu)$ is a coboundary quasi-differential algebra.
In the adjoint representation the zero-degree element $I$
corresponds to the grading character:
$$ad_I(f)=-\partial eg(f) f$$
and Bianchi's identity holds:
$$d\alpha_\mu=0, \qquad \alpha_\mu(x,y,z)=(xy)z-x(yz)$$
Moreover
$$d^2 s=[\alpha,s]$$
\end{prop}
\begin{pf}
Note first that the identity map $I:A\to A$ has degree zero. If
$f\in C^p(A)$
\begin{align}
ad_I(f) &=I\circ f -(-1)^{p\cdot 0} f\circ I\notag \\
        &=f - (p+1) f\notag\\
        &=- \partial eg(f) f
\end{align}
Since $[I,f]=-[f,I]$, the right adjoint action of the unit is scalar
multiplication by the degree map.

Obvoiusly $d_\mu I=[\mu, I]=\partial eg(\mu) \mu=\mu$, thus I is a unit.
Now, by the Jacobi identity, the associator is a cocycle
$d_\mu \alpha=0$ (the curvature is closed):
$$ d_\mu [\mu,\mu]=[\mu,[\mu,\mu]]=0 $$
where the assumtion that 2 and 3 do not annihilate non-zero elements in A
was used.

The second equation follows from lemma \ref{L:1}.
\end{pf}
In the context of the previous proposition, for any element $x$
of even degree, we have  $[x,x]=0$, since $0=[I,[x,x]]=-2 \partial eg(x) [x,x]$
(see lemma \ref{L:1}).

To define the cohomology of certain classes of non-associative algebras,
we give the following definition:
\begin{defin}
An algebra $(A,\mu)$ is called {\em N-coherent} if $d_\mu^N=0$.
\end{defin}
Note that an algebra is associative iff it is a 2-coherent algebra,
and a 1-coherent algebra is just the trivial one:$\mu=0$.

Of course an algebra $(A,\mu)$ is N-coherent iff $ad_\mu$ is a nilpotent
element of order $N$.

For an N-coherent algebra the quasi-differential graded Lie algebra
$(C^\bullet,d_\mu)$ is an N-complex as defined in \cite{Ka}.

%
%
Next we will consider examples of non-associative algebras.
\begin{example}\label{E:Liealg}
If $(A,\mu=[,])$ is a Lie algebra, then its associator is:
\begin{align}
\alpha(a,b,c) &=([[a,b],c]-[a,[b,c]])\notag\\
              &=[b,[c,a]]
\end{align}
and $Alt(\alpha)=0$ as expected.
\end{example}
\begin{example}
If $(A,m)$ is an associative algebra, we can associate
to it the Jordan algebra $(A,\mu_+)$, with
$\mu_+(a,b)=\{a,b\}=ab+ba$,
and its Lie algebra $(A,\mu_-)$, with
$\mu_-(a,b)=[a,b]=ab-ba$.
The corresponding associators are
$$\alpha_+(a,b,c)=\{cab+bac-acb-bca\}=-[b,[c,a]]$$
and
$$\alpha_-(a,b,c)=\{-bac-cab+acb+bca\}=[b,[c,a]]$$
Thus they have the same even quasi-differential Lie algebra
$(C^+,d_{\alpha_\pm})$ (modulo a sign).

Note also that $Alt(\alpha_\pm)=0$.
\end{example}
%

\section{Torsion Algebras}\label{S:TA}
Hochschild differential complex
is defined for an associative algebra with
coefficients in a symmetric $A$-bimodule $M$.
When relaxing both conditions,
associativity and action requirement, one obtains formulas
which are familiar in differential geometry, 
and correspond to a non-flat connection.
It may be thought of as a geometry of ``vector fields''
without starting from a function algebra. 
\begin{defin}\label{D:ta}
A {\em torsion algebra} $\C{M}=(C,\mu,[,]_C))$ is a (non-associative)
k-algebra $(C,\mu)$ together with a Lie bracket $[,]_C$.
Its {\em torsion} is $T=\mu_--[,]_C$, where $\mu=\mu_++\mu_-$ is the
decomposition into its {\em symmetric (quasi-Jordan)} and
{\em skewsymmetric (quasi-Lie)} part:
$$T(X,Y)=\mu(X,Y)-\mu(Y,X)-[X,Y]_C, \qquad X, Y \in C.$$
\end{defin}
It is a generalization of the most important classes of algebras.
Associative algebras, with the usual Lie bracket $[x,y]=xy-yx$ are
torsion algebras, with $T=0$. Lie algebras $(C,[,])$, with
$\mu=\frac12[,]$ are again torsion algebras with zero torsion.
Poisson algebras (compatibility between $\mu$ and $[,]$)
and Gerstenhaber algebras (non-commutative Poisson algebras) 
can be interpreted as torsion algebras in several ways.
Pre-Lie algebras (non-associative algebras such that
the skew part of the multiplication is a Lie bracket) with
$[,]=\mu_-$ are torsion algebras with zero torsion ($T=0$).

We think of $(C,[,])$ as a Lie algebra of vector fields and the
multiplication $\mu$ as a {\em generalized connection}.
\begin{example}\label{Ex:1}
Obviously any manifold $V$ with a connection $D$ defines a torsion algebra.
Take $C$ as the Lie algebra of vector fields on $V$ and interprete
the connection as a non-associative multiplication
$\mu(X,Y)=D_XY$.
Then the torsion is $T=D_--[,]$, i.e.
$$T(X,Y)=D_XY-D_YX-[X,Y].$$
In this {\em geometric example}, the torsion tensor
coincides with the torsion in the sense of definition \ref{D:ta}.

If $V$ is the real line, then the Lie algebra of vector fields
$X_f=f\partial_t$ can be
identified with $(C^\infty(V),[,])$, where $[f,g]=fg'-gf'$.
Also any connection $D$ has a canonical Christoffel symbol $\Gamma$ and
\begin{equation}\label{Eq:1}
D_fg=f(g'+g\Gamma)\qquad (D=d+\Gamma).
\end{equation}
\end{example}

\subsection{Notation and background}
$(C,\mu)$ will denote a possibly non-associative $k$\@-algebra, where
$k$ is a ring. We will write $D_XY=\mu(X,Y)$, in order
to emphasize the geometric interpretation.
Basic definitions for the usual algebraic model are assumed,
following \cite{Ko}.
The prefix $\C{M}$ will be used with notions referring to the
formal context (``non-commutative'' space),
and the prefix $\C{A}$ to refer to the usual notions in the context
of a geometric example, e.g. on a manifold $V$.

In the ``geometric world'', functions can be identified as
$k$\@-endo\-mor\-phisms (multiplication of vector fields by functions) for
which the connection is linear in the first argument.
%
%
\subsection{The Algebra of Functions}
We define as ``functions'' the annihilator of the left commutator
of the multiplication $D$:
\begin{defin}
Let $(C, D, [,]_C)$ a torsion algebra.
Its elements are called {\em $\C{M}$-vector fields}.
The set of {\em $\C{M}$-functions} is:
$$A=\{\phi\in End_k(C)| D_{\phi(X)}Y=\phi(D_XY)\}$$
The {\em multiplication} of $\C{M}$-functions is the natural composition
of $k$\@-endo\-mor\-phisms in $End_k(C)$.
\end{defin}
Note that the multiplication of $\C{M}$-functions is an internal operation:
$$D_{(\phi\circ\psi)x}=D_{\phi(\psi x)}=\phi D_{\psi x}=\phi \psi D_x$$
and that the set of $\C{M}$-vector fields $C$ is a left $A$-module.

The multiplication $D$ defines a $k$-linear map:
$$\tilde{D}:C\to End_k(C)$$
called the {\em left regular quasi-representation} of $(C,D)$
as a non-associative algebra.

We will test the notions introduced against the simple geometric
example of the real line.
\begin{example}
In the context of example \ref{Ex:1}, multiplication of vector fields
$C\cong C^\infty(V)$ by functions is just the regular left representation
$L:C^\infty(V)\to End(C^\infty(V))$ of $C^\infty(V)$ (in the usual sense):
$$(fX_g)=f(g\partial_t)=(fg)\partial_t=X_{fg}$$
Moreover the {\em $\C{M}$-functions} $A$ are naturally identified as
{\em $\C{A}$-functions} $C^\infty(V)$.
Indeed, if $\phi\in End_k(C)$ ``left commutes'' with $D$:
$$D_{\phi(f)}g=\phi(D_fg)$$
then (see equation \ref{Eq:1}):
$$ \phi(f)(g'+g\Gamma)=\phi(f(g'+g\Gamma)).$$
But it is clear that $g'+g\Gamma=h$ has a solution for any
$h\in C^\infty(V)$.
Thus $\phi(fh)=\phi(f)h$, so $\phi(h)=\phi(1)h$ and $\phi$ 
corresponds to left multiplication by $\phi(1)$.

We note that $\phi$ is a function iff
$\tilde{D}\circ\phi=L_\phi\circ \tilde{D}$, where:
$$L:End(C)\to End(End(C))$$
is the regular representation of the associative algebra $(End(C),\circ)$.
In other words, $\tilde{D}$ intertwines $\phi$ and $L_\phi$:
$$\tilde{D}\circ \phi=L_\phi\circ \tilde{D}$$ 
\end{example}
We will interpret the $\C{M}$-vector fields as derivations on $A$.
Let $X\in C$ and $\phi\in A$. 
\begin{lem}\label{L:21}
Any two of the following conditions imply the third:\\
(i) The action of $C$ on functions is defined by:
\begin{equation}
(X\cdot\phi)(Y)=[X,\phi(Y)]_C-\phi([X,Y]_C), \qquad Y\in C
\end{equation}
(ii) $D$ is a derivation law:
\begin{equation}
D_X(\phi Y)=(X\cdot\phi)Y+\phi D_XY, \qquad ( X\cdot\phi=[D_X,\phi])
\end{equation}
(iii) The torsion is $A$-bilinear.
\end{lem}
\begin{pf}
Note that the torsion $T$ is skewsymmetric and:
$$\begin{array}{rll}
T(X,\phi Y)&=&D_X(\phi Y)-D_{\phi Y}X -[X,\phi Y]\\
&=&\{D_X(\phi Y)-\phi D_XY -(X\cdot \phi)Y\} \
  +\ \phi T(X,Y)\\
& &+\ \{(X\cdot\phi)Y+\phi[X,Y]-[X,\phi Y]\}
\end{array}$$
Now it is clear that any two conditions imply the third:
$$\begin{array}{rl}
T(X,\phi Y)-\phi T(X,Y)=& \{D_X(\phi Y)-\phi D_XY -(X\cdot \phi)Y\}\\
 & +\ \{(X\cdot\phi)Y+\phi[X,Y]-[X,\phi Y]\}
\end{array}$$
\end{pf}
We will adopt the second condition in lemma \ref{L:21} as a definition
for the action of a vector field on a function.
\begin{defin}\label{defvf}
A vector field $X\in C$ acts on a function $\phi\in A$ by:
\begin{equation}\label{eqvf}
X\cdot\phi=[D_X,\phi]
\end{equation}
\end{defin}
Note that it measures the failure of $D$ to be right $A$-linear.
\begin{prop}
The $\C{M}$-vector fields act as (external) derivations on $A$.
\end{prop}
\begin{pf}
If $\phi$ and $\psi$ are $\C{M}$-functions, then:
$$(X\cdot(\phi\circ\psi))Y=D_X(\phi(\psi(Y)))-(\phi\circ\psi)D_XY$$
and
$$ (X\cdot\phi)\circ\psi (Y)+\phi\circ(X\cdot\psi)(Y)=
D_X(\phi(\psi(Y)))-\phi(\psi(D_XY)) $$
\end{pf}
For $X\cdot\phi$ to be again a function,
so that elements of $C$ act as derivations, we note the following:
\begin{lem}\label{L:22}
The following conditions are equivalent:

(i) For any $X\in C$ and $\phi\in A$, $X\cdot\phi$ is an $\C{M}$-function.

(ii) The associator $\alpha$ of $D$ is $A$-linear in the first two
variables.
\end{lem}
\begin{pf}
Recall that the associator is:
$$\alpha(X,Y,Z)=(XY)Z-X(YZ)=D_{D_XY}Z-D_XD_YZ$$
where multiplicative notation was alternatively used.
The following are equivalent:
$$D_{(X\cdot\phi)Y}Z=(X\cdot\phi)D_YZ$$
$$D_{D_X\phi Y}Z-\phi D_{D_XY}Z=D_X(\phi D_YZ)-\phi D_XD_YZ$$
$$(X\circ\phi(Y))\circ Z-X\phi(Y\circ Z)=\phi(\alpha(X,Y,Z))$$
$$(X\circ\phi(Y))\circ Z-X\circ(\phi(Y)\circ Z)
+X\circ[\phi(Y)\circ Z -\phi(Y\circ Z)]=\phi(\alpha(X,Y,Z))$$
$$\alpha(X,\phi(Y),Z)+X\circ[D_{\phi(Y)}Z-\phi(D_YZ)]=\phi(\alpha(X,Y,Z)).$$
Since:
$$D_{D_{\phi X}Y}Z=D_{\phi D_XY}Z=\phi (D_{D_XY}Z)$$
the linearity of the associator in the first variable is clear:
$$\begin{array}{rl}
\alpha(\phi X,Y,Z)&=D_{D_{\phi X}Y}Z-\phi(D_XD_YZ)\\
&=\phi \alpha(X,Y,Z).
\end{array}$$
A direct computation proves the $A$-linearity in the second variable:
$$\begin{array}{rl}
\alpha(X,\phi Y,Z)&=D_{D_X\phi Y }Z-D_XD_{\phi Y}Z
=D_{D_X\phi Y}Z-D_X(\phi D_YZ)\\
&=D_{((X\cdot\phi)Y+\phi D_XY)}Z-((X\cdot\phi)D_YZ+\phi D_XD_YZ)\\
&=\phi(D_{D_XY}Z-D_XD_YZ)=\phi\alpha(X,Y,Z)
\end{array}$$
\end{pf}
\begin{defin}
A torsion algebra $\C{M}=(C,D,[,]_C)$ is called {\em regular} if
the torsion and the associator are $A$-bilinear in the first two
variables.
\end{defin}
Since condition (ii) of lemma \ref{L:21} holds by definition
in a torsion algebra, any of the other two imply the third.
Now we can easily prove the following:
\begin{theorem}\label{th1}
Let $\C{M}=(C,D,[,]_C)$ be a regular torsion algebra.
Then, for all $X,Y\in C$ and $\phi\in End_k(C)$ an $\C{M}$-function:

(1) $X\cdot\phi$ is a $\C{M}$-function.

(2) $X$ acts as a derivation on $\C{M}$-functions.

(3) The associator of $D$ is $A$-linear in the first two variables.

(4) $[X,\phi Y]=\phi[X,Y]+(X\cdot\phi)Y$

(5) $D$ is a connection on $\C{M}$: $D_X(\phi Y)=(X\cdot\phi)Y+\phi D_XY$.

(6) The torsion $T$ is $A$-bilinear.
\end{theorem}
\begin{pf} 
(1) holds by assumption.

Since by definition the conditions (ii) and (iii) from lemma \ref{L:21}
hold, (2) follows.

The other statements are clear from lemma \ref{L:21} and \ref{L:22}.
\end{pf}
In what follows we will omit the $\C{A}/\C{M}$ prefix.
\subsection{Differential Forms}
The exterior derivative will be defined as the differential
of the Cheva\-lley-Eilenberg quasi-complex.

In the context of non-associativity, it is natural to relax the action
requirement as well.
The associator (algebraic point of view) may be
interpreted as a curvature (geometric point of view) or as a
monoidal structure through categorification (failure to be a morphism)
(see \ref{S:I}).
Now, an action of $A$ on $M$ is a morphism $\rho:A\to End(M)$ and
an associative multiplication is an action $A\to End(A)$.
\begin{defin}
A {\em quasi-action} of $A$ on $M$, in the category of $k$-modules,
is a $k$-linear map $L:A\to End_k(M)$.
\end{defin}
Let $M$ be a left $A$-module with a derivation law $D^M:C\to End_k(M)$.
\begin{defin}
The {\em $M$-valued $\C{M}$-differential forms} are defined as usual:
$$\Omega^n(\C{M},M)=\{\omega:C\times...\times C\to M |\ \omega\
\text{alternating and\ } A-\text{multilinear}\}$$
\end{defin}
Then $\Omega(\C{M},M)$ is just the alternate part of the Hochschild
cochains $C^\bullet(C;M)$ with coefficients in $M$ (Chevalley cochains).

To define first the Hochschild quasi-complex ($d^2$ not necessarily zero),
consider the following $C$-quasi-bimodule structure on $M$:
$$\diagram
 \lambda:C\times M\to M, & \lambda(X,u)=D^M_Xu &
   C\times M \dto_{-\sigma_{(12)}} \drto^{\lambda=D^M} & \\
 \rho:M\times C\to M, & \quad\rho(u,X)=-D^M_Xu &
   M\times C \rto_{\rho=\lambda^{op}} & C
\enddiagram$$
where $\lambda^{op}$ is the opposite quasi-action using the signed braiding.
In the associative case with $M=A$, the use of the signed
braiding gives $M$ a structure of
$(A,A_{op})$ {\em supersymmetric bimodule} structure: $am=-ma$.

Instead of the Hochschild quasi-complex
derived from the associated graded Lie algebra $(C^\bullet(C),[,])$,
with $d\omega=[\mu, \omega]=\mu\circ \omega -(-1)^p\omega\circ \mu$,
consider the Hochschild quasi-complex
$C^p(C;M)=Hom_A(C^p, M)$, $p\ge-1$ of the
Lie algebra $(C,[,]_C)$ as a non-associative algebra,
with coefficients the $C$-quasi-bimodule $M$:
$$d\omega=(-1)^p((\lambda,\rho)\circ \omega -(-1)^p \omega\circ [,]_C)$$
$$\begin{array}{rl}
d\omega(a_1,...,a_{p+2})=&
 \lambda(a_1,\omega(a_2,...,a_{p+2}))-\omega([a_1,a_2]_C,...,a_{p+2})+...\\
 & +(-1)^p \rho(\omega(a_1,...,a_{p+1}),a_{p+2})
\end{array}$$
Then for $u\in C^{-1}=M$ and $\omega\in C^0$:
$$\begin{array}{rl}
du(X)&=\lambda(X,u)-\rho(u,X)=2D_Xu\\
d\omega(X,Y)&=D_X \omega(Y)-D_Y\omega(X)-\omega([X,Y]_C)\\
ddu(X,Y)&=D_X du(Y)-D_Y du(X)-du([,]_C)\\
&=2 (D_X D_Yu-D_Y D_Xu-D_{[X,Y]_C})u=K(X,Y)u
\end{array}$$
To obtain the usual formulas in geometry, consider the alternating
part $\Lambda^\bullet(A;M)$ of the above complex,
and project the differential $d_{Ch}=Alt\ \circ\ d$.
\begin{defin}
$(\Lambda^\bullet(A;M), d_{Ch})$ is called the associated
Chevalley-Eilenberg quasi-complex of $C$ with coefficients in $M$.
\end{defin}
Then, for example:
$$d_{Ch}\omega(X,Y,Z)=\sum_{cycl} D_X\omega(Y,Z)-\omega([X,Y],Z)$$
\subsection{The Lie Derivative}
Let $\C{M}=(C,D,[,]_C)$ be a regular torsion algebra. 
Consider the $A$-module $M=A$ and the corresponding differential
forms $\Omega(\C{M})$. The {\em canonical derivation law} on $A$ is:
$$D_X\phi=X\cdot\phi.$$
Extend as usual the Lie derivative defined on functions
and vector fields as a derivation on the tensor algebra
commuting with contractions.
It is easy to see that it is an internal operation.
For example, if $\omega:C\to A$ is a 1-form, then
$(\C{L}_X\omega)(Z)=D_X\omega(Z)-\omega([X,Z])$ is $A$-linear.

An exterior differential on forms $\Omega(A;\C{F})$ is defined by
the homotopy formula: $\C{L}_X=d i_X+i_X d$.
The usual explicit formula holds for $d$.
It coincides with $d_{Ch}$ defined above.


%
%
\section{Examples}\label{S:Ex}
We will start by investigating torsion algebras for which $T=0$.

Let $(C, D)$ be a unital associative algebra.
Consider the corresponding Lie algebra structure: $[X,Y]_C=D_XY-D_YX$.
Then the torsion is $T=0$.
The associator is zero and $(C,D,[,])$ is a regular torsion algebra.
If $\phi\in End_k(C)$ is a function, then $D_{\phi(X)}Y=\phi(D_XY)$
in multiplicative notation is just $\phi(X)Y=\phi(XY)$. Thus
$\C{M}$-functions are left multiplication by elements of $C$ and
the algebra of functions is isomorphic to the initial algebra.
The morphism $C\to Der(A)$, realizing $C$ as derivations of $A$,
is the usual Lie algebra representation.

Thus we have the following:
\begin{prop}
Any associative algebra $(C,\mu)$ has a natural structure of a
regular torsion algebra.
The algebra $C$ is isomorphic with the algebra of $\C{M}$-functions.
\end{prop}

\newcommand{\cbar}{\overline{\circ}}
\newcommand{\ctld}{\tilde{\circ}}

\subsection{Hochschild Pre-Lie Algebra}
Let $V$ be a $k$-module and $C=(C^\bullet(V), \cbar)$ the corresponding
graded pre-Lie algebra, where $\cbar$ denotes the Gerstenhaber composition
operation \cite{Ge}.
Recall that $[x,y]=x\cbar y-(-1)^{pq}y\cbar x$, defines a
Lie bracket, where $x\in C^p(V)=Hom_k(V^{p+1},V)$ and
$y\in C^q(V)=Hom_k(V^{q+1},V)$ are two Hochschild cochains.

Then $C$ is a torsion algebra with $D=\cbar$ and $T=D_--[,]=0$.

A $k$\@-endomorphism $\phi\in End_k(C)$ is an $\C{M}$-function iff:
\begin{equation}\label{Eq:2}
D_{\phi x}y=\phi(D_xy)
\end{equation}
and an argument similar to the case of associative algebras gives
$\phi=L_{\phi(1)}$, where $1=id_V$ and
$L:C\to (End_k(C),\ctld)$ is the regular
quasi-representation.
Note that $\cbar$ is not associative and $L_{id_V}$ is only
a projector on the even part of $C$.

Denote $\phi(1)$ by $f$. Then equation \ref{Eq:2} holds iff:
$$ D_{f\cbar x}y=f\cbar (D_xy)$$
i.e. $(f\cbar x)\cbar y=f\cbar(x\cbar y)$ for any $x,y\in C^\bullet$.
It is easy to see that this is true iff $f\in C^0(V)$ and thus
the set of functions is $A=C^0(V)$.

The composition of functions is composition of $k$-endomorphisms:
$$L_f\ctld L_g=L_{f\cbar g}, \qquad f,g\in A.$$
since $\cbar$ reduces to the usual composition $\circ$ of
$k$\@-endomorphisms of $C^0(V)=End_k(V)$.
Thus we have:
\begin{theorem}\label{th2}
Let $V$ be a $k$-module and $(C^\bullet(V),\cbar)$ the corresponding
pre-Lie algebra. Then:

(i) $C=(C^\bullet(V), D, [,])$ is a zero torsion algebra where
$D=\cbar$ is called the canonical connection.

(ii) The algebra of functions of $C$ is $A=(C^0(V),\cbar)$, i.e.
$(End_k(V),\circ)$.

(iii) $C$ acts through exterior derivations on $A$:
$$\begin{array}{rl}
(x\cdot L_f)(y)&=D_x(L_f(y))-L_f(D_xy), \quad x, y, f \in C
\end{array}$$
where $L:C\to End_k(C)$ is the regular left quasi-representation of $C$.
\end{theorem}
We note that the failure to be a regular torsion algebra
comes from the $A$ non-linearity of the associator.
Recall that the associator is graded skewsymmetric in the
last two variables (\ref{L:comp1}). Thus being a regular torsion algebra
would be equivalent to $\alpha$ being an $A$-multilinear form.
\subsection{Poisson Algebras}
Let $(C,\cdot,\{,\})$ be a Poisson algebra, with $D=\cdot$ commutative and
associative and Lie bracket $[,]$ being the Poisson bracket $\{,\}$.
Then $(C, D, [,])$ is a torsion algebra with torsion $T=-[,]$.
Since $D$ is associative, its algebra of $\C{M}$-functions $A$ is
isomorphic to $C$, in a manner similar to the associative algebra case.
In this way, a Poisson algebra is not a regular torsion algebra.

If $D=[,]=\{,\}$ then it becomes a zero torsion algebra, but it is not
clear what the algebra $(A,\circ)$ of $\C{M}$-functions is and what is
the relation with the multiplication of functions.



\begin{thebibliography}{11}

\bibitem[DM]{DM} Michel Dubois-Violette, Thierry Masson,
``$SU(n)$-connections and noncommutative differential geometry,
{\em dg-ga/9612017}.

\bibitem[DMM]{DMM} M. Dubois-Violette, J. Madore, T. Mason, J. Mourad,
``On curvature in Non-commutative Geometry'',
{\em q-alg/9512004}.

\bibitem[Ge]{Ge} Gerstenhaber, M., Schack, S. D.:
Algebraic Cohomology and Deformation Theory,
{\em Deformation theory of algebras and structures and applications},
(Il Ciocco, 1986), 11-264,
NATO Adv. Sci. Inst. Ser. C: Math. Phys. Sci., 247,
Kluwer Acad. Publ., Dordrecht, 1988.

\bibitem[GS]{GS} Gerstenhaber, M., Schack, S.D., ``Algebras, Bialgebras,
Quantum Groups, and Algebraic Deformations'', Contemporary Mathematics,
Vol. 134, 1992, 51-92.

\bibitem[I]{I}
Ionescu, Lucian M., ``On Categorification'',
{\em math.CT/9906038}, (1999).

\bibitem[Ka]{Ka} Kapranov, M.M., ``On the q-analog of Homological Algebra'',
q-alg/9611005.

\bibitem[Ko]{Ko} Koszul, J.,L., ``Lectures on Fibre Bundles and
Differential Geometry'', Springer-Verlag, Berlin, 1986.

\bibitem[KPS]{KPS} Liivi Kluge, Eugen Paal, Jim Stasheff,
``Invitation to composition'',
{\em math.QA/9909083 v2}.

\bibitem[W]{W} Weibel, C.A., ``An Introduction to Homological Algebra'',
Cambridge University Press, 1994.

\end{thebibliography}
\end{document}